\documentclass[english,12pt,a4paper]{article}
\usepackage{amsmath,amssymb,amsthm} 

\usepackage{graphicx}
\usepackage[scanall]{psfrag}
\usepackage{comment}
\usepackage{tipa}

\newcommand{\bA}{\boldsymbol{A}}
\newcommand{\bB}{\boldsymbol{B}}

\newcommand{\by}{\boldsymbol{y}}

\newcommand{\abs}[1]{\left | {#1} \right |}

\title{How far are vowel formants from \\ computed vocal tract resonances?}

\author{D.~Aalto$^{1,3}$,  A.~Huhtala$^2$,  A.~Kivel\"a$^2$, J.~Malinen$^2$, \\  P.~Palo$^3$, J.~Saunavaara$^4$, M.~Vainio$^3$
\\ {\small $^1$ Dept. of Signal Processing and Acoustics, Aalto
  University, Finland } \\ {\small $^2$ Dept. of Mathematics and Systems
  Analysis, Aalto University, Finland } \\ { \small $^3$ Inst. of Behavioural
  Sciences (SigMe group), University of Helsinki, Finland } \\ { \small $^4$
  Medical Imaging Centre of Southwest Finland }  \\ {\small \tt
    jarmo.malinen@aalto.fi }}

\date{\today}

\begin{document}

\maketitle

\begin{abstract}
We compare numerically computed resonances of the human vocal tract
with formants that have been extracted from speech during vowel
pronunciation. The geometry of the vocal tract has been obtained by
MRI from a male subject, and the corresponding speech has been
recorded simultaneously. The resonances are computed by solving the
Helmholtz partial differential equation with the Finite Element Method
(FEM). 

Despite a rudimentary exterior space acoustics model, i.e., the
Dirichlet boundary condition at the mouth opening, the computed
resonance structure differs from the measured formant structure by
$\approx$ 0.7 semitones for \textipa{[i]} and \textipa{[u]} having
small mouth opening area, and by $\approx$ 3 semitones for vowels
\textipa{[a]} and \textipa{[ae:]} that have a larger mouth
opening. The contribution of the possibly open velar port has not been
taken into consideration at all which adds discrepancy for
\textipa{[a]} in the present data set.  We conclude that by improving
the exterior space model and properly treating the velar port opening,
it is possible to computationally attain four lowest vowel formants
with an error less than a semitone. The corresponding wave equation
model on MRI-produced vocal tract geometries is expected to have a
comparable accuracy.
\end{abstract}

\noindent {\bf Keywords.} Formant analysis, acoustic resonance computation, FEM, MRI.

\section{Introduction}

The purpose of this paper is to evaluate the accuracy of vowel
simulations based on the wave equation model \eqref{eq:model}.
We use 3D vocal tract (VT) geometries that have been obtained by
Magnetic Resonance Imaging (MRI) from a native male speaker of Finnish
while he pronounces prolonged vowels \textipa{[\textscripta]},
\textipa{[i]}, \textipa{[u]}, and \textipa{[\oe]}.
The evaluation is carried out by comparing the computed resonances of
\eqref{eq:model} with the measured formants, extracted from sound
samples, instead of comparing simulated and actual speech signals.  In
this work, \emph{the sound samples have been recorded simultaneously
  with the MRI}, using the equipment and the arrangements detailed in
\cite{A-H-M-P-P-S-V:RSSAMRI,A-M-V-S-P:EMAEMRID,M-P:RSDMRIII,Palo:L:2011}.
This is in contrast to, e.g., our earlier work \cite{H-L-M-P:WFWE}
where only a single anatomic configuration (corresponding to Swedish
\textipa{[\o]}) was taken from the data set of
\cite{Engwall:CA23D:1999}.

We use the same wave equation model for vowels as in
\cite{H-L-M-P:WFWE}, namely
\begin{equation}\label{eq:model}
  \begin{cases}
     & \Phi_{tt}  =  c^2\Delta \Phi    
    \text{ on }  \Omega,  
    \Phi=0   \text{ on }  \Gamma_1, \\ & 
    \frac{\partial \Phi}{\partial\nu}=0 
    \text{ on }  \Gamma_2,    
    \Phi_t + c \frac{\partial \Phi}{\partial\nu}= 2 \sqrt{\frac{c}{\rho_0}} u  \text{ on } \Gamma_3,
  \end{cases}
\end{equation}
where $\Omega \subset \mathbb R^3$ is the interior of the VT whose
boundary $\partial \Omega = \overline{\Gamma_1} \cup
\overline{\Gamma_2} \cup \overline{\Gamma_3}$ consists of the mouth
opening $\Gamma_1$, the VT tissue walls $\Gamma_2$, and an (imaginary)
control surface $\Gamma_3$ placed right above the glottis. The
parameters $c = 350$ m/s  and $\rho_0 = 1.225$ kg/m$^3$ are the speed of
sound in and the density of dry air at $T = 305$ K, respectively. The
functions in \eqref{eq:model} are as follows: $u = u({\bf r},t)$ is
the incoming acoustic power (per unit area) at the glottis input,
$\frac{\partial \Phi}{\partial\nu} = {\bf \nu} \cdot \nabla \Phi$, and
${\bf \nu}$ is the exterior unit normal on $\partial \Omega$.  In time
domain simulations, we compute the velocity potential $\Phi({\bf r},t)$
for a given glottal input function $u({\bf r},t)$ produced by a source
such as described in \cite{A-A-M-V:MLBVFVTO,A-A-M:LFPSGFM}.
 From the velocity potential, the sound pressure and (perturbation)
 velocity can be extracted as the partial derivatives
$p' = \rho_0 \Phi_t$  and ${\bf v} = \nabla \Phi$.
The physics of model \eqref{eq:model} is further explained in
\cite{H-L-M-P:WFWE} and the references therein.

In the past, the VT acoustics has been modelled in many different
ways. Electrical \emph{transmission lines} have been used already in
\cite{Dunn:CVR:1950}, and the classical \emph{Kelly--Lochbaum model}
in \cite{Kelly:SS:1962} makes use of reflection/transmission
coefficients of a variable diameter tube. The 3D \emph{wave equation}
for linear wave propagation as well as the related \emph{Helmholtz
  equation} for acoustic resonances have been known for a long time;
see, e.g. \cite{Helmholtz:LT:1863}. The Kelly--Lochbaum model is
closely related to \emph{Webster's equation} in, e.g.,
\cite{Chiba:VNS:1958,Fant:ATS:1960} but the latter can be easily
deduced from variable-impedance electrical transmission lines as well
as from the wave equation in 3D tubular domains as shown in great
generality in \cite{L-M:WECAD}.  More advanced models are the
\emph{transmission line networks} that have been applied for speech
in, e.g., \cite{ElMasri:VTA:1996, ElMasri:DTL:1998, Mullen:WPM:2006};
see also \cite{A-M:CPBCS} for a purely mathematical treatment.

At their best, all of these modelling paradigms are known to produce
very good simulated speech even though they are based on radically
simplified representations of the underlying anatomic geometry
$\Omega$ with the sole exception of the wave equation.  In most
applications, simplifications are even desirable as it may improve
conceptual clarity and reduce the computational burden. There are,
however, situations where the faithful representation of the VT
geometry is required, e.g., when modelling the effects of anatomical
abnormalities and maxillofacial surgery on speech
\cite{Dedouch:FEM:2002,Niemi:AC:2006,Nishimoto:ETF:2004,T-M-K:AAVTDVPFDTDM,Vahatalo:EG:2005,Svancara:NME:2006}.

\section{Background and motivation}

In our earlier work \cite{H-L-M-P:WFWE}, the same numerical
computations were carried out using a minimal data set: a single
MRI-based anatomic geometry $\Omega$ corresponding to the Swedish
vowel \textipa{[\o]}.  The computed resonances were compared to the
first four formants of all Swedish vowels (including \textipa{[\o]})
that were extracted from speech samples of the same test subject. The
speech samples were not recorded simultaneously with the MRI because
of technological restrictions but the subject was in a similar supine
position during both MRI and speech recording; see
\cite{Engwall:CA23D:1999}. 

We made the following observations in \cite{H-L-M-P:WFWE}:
\begin{enumerate}
\item The computed resonances $R_1$...$R_4$ corresponding the formants $F_1$...$F_4$
  of \textipa{[\o]} are systematically $3\tfrac{1}{2}$ semitones too
  high compared to the measured values;
\item the \emph{formant ratios} of the computed and measured data
  (i.e., $R_i/R_1$ and $F_i/F_1$ for $i = 2,3,4$) correspond to each
  other quite well; and
\item if the systematic error in $R_1$...$R_4$ of \textipa{[\o]}
  is compensated by linear scaling, then the scaled, computed 
  data gets identified correctly as \textipa{[\o]} in the measured
  formant table from the same subject.
\end{enumerate}


Two main potential sources were identified for the discrepancy: (i)
the Dirichlet boundary condition on $\Gamma_1$ in \eqref{eq:model}
results in \emph{too short acoustic length of the computational VT};
and (ii) the minimal \emph{data set used in \cite{H-L-M-P:WFWE} is
  insufficient} to draw any conclusions on the error sources.  The
purpose of this work is to exclude the latter possibility (ii) by
extending and improving the data set in an essential manner. We also
aim at a deeper understanding of the sources of descrepancy to guide
future model improvements and to understand the quality of simulation
that one can reasonably expect to attain.

We remark that the formant computation of \cite{H-L-M-P:WFWE} was
later validated by independent FEM resonance computations that were
based on the generalized Webster's model instead of \eqref{eq:model};
the computed resonances $R_1$...$R_4$ corresponding to $F_1$...$F_4$
were practically the same as reported in \cite[Table 3.1 on
  p. 31]{Aalto:dt:2009}. As shown in \cite{L-M:WECAD}, the generalized
Webster's model is a low-frequency approximation of the wave equation
in a tubular domain. In the case of human VT, the approximation
remains accurate for formants $F_1$...$F_3$ and even for $F_4$ at
least in some vowel configurations where cross-mode resonances do not
dominate; see \cite[Fig. 1]{H-L-M-P:WFWE}.

We comment on the interesting parallel work \cite{T-M-K:AAVTDVPFDTDM}
at the end of the article.

\section{Model and methods}

As explained in \cite[p. 3]{H-L-M-P:WFWE}, the resonances of
Eq. \eqref{eq:model} can be solved by finding the complex frequencies
$\lambda$  such that the Helmholtz problem
\begin{equation}
  \label{eq:eigenfuncs}
  \begin{cases}
    \lambda^2 \Phi_\lambda=c^2\Delta \Phi_\lambda \text{ on } \Omega, 
    \Phi_\lambda = 0 \text{ on } \Gamma_1, \\ 
    \frac{\partial \Phi_\lambda}{\partial\nu}=0 \text{ on } 
    \Gamma_2, \text{ and } 
    \lambda \Phi_\lambda + c\frac{\partial \Phi_\lambda}{\partial\nu}= 0 
    \text{ on } \Gamma_3
\end{cases}
\end{equation}
is solvable for some nonzero eigenfunctions $\Phi_\lambda({\bf
  r})$. It is known that all such eigenvalues $\lambda$ form an
infinite sequence $\{\lambda_j \}_{j \in \mathbb Z\setminus \{ 0 \}}$
with $\abs{\lambda_j} \to \infty$ as $\abs{j} \to \infty$,
$\mathop{Re}{\lambda_j} \leq 0$, and $\mathop{Im}{\lambda_{-j}} = -
  \mathop{Im}{\lambda_{j}}$. The lowest formants $F_1$, $F_2$,
    $\ldots$, correspond to the numbers $R_j = \mathop{Im}{\lambda_j}$ for
    $j = 1, 2, \ldots$ in the order of increasing imaginary parts.

As solving Eqs.~\eqref{eq:model} and \eqref{eq:eigenfuncs}
analytically is possible only in a radically simplified geometry
\cite{Sondhi:RBV:1986}, we solved the problem numerically by the
\emph{Finite Element Method} (FEM).  This is the approach used by
\cite{Lu:FES:1993},
\cite{Niikawa:FDV:2002}, \cite{Dedouch:AMA:2002},
\cite{Sasaki:FEM:2003}, \cite{Svancara:NMP:2004}, and by many others.
Eqs. \eqref{eq:eigenfuncs} were solved in variational form as given
in \cite[Eq.~(5)]{H-L-M-P:WFWE} using a custom implementation of FEM
programmed in MATLAB. We used piecewise linear shape functions on
tetrahedral meshes. The tetrahedral meshes were generated using TetGen
\cite{TetGen:2011} from a triangular surface mesh. As a result, we
obtained the matrices $\bA$ and $\bB$ for a high-order eigenvalue
problem $\bA \by(\lambda) = \lambda \bB \by(\lambda)$ as explained in
\cite[Eq.~(6)]{H-L-M-P:WFWE}. The lowest eigenvalues $\lambda_j$,
$j = 1,2,3,4$ were then computed using the \emph{eigs} routine of MATLAB.
It takes around 3 seconds on a workstation with an Intel Xeon X3450
processor to build the matrices and to solve the eigenvalue problem.

The imaginary parts of the computed $\lambda_j$ are given in
Table~\ref{numberofelements} together with the number of elements
used in each VT geometry. The computed values are good approximations
of the eigenvalues defined in Eqs.~\eqref{eq:eigenfuncs} when the number
of elements is high enough. It was observed with the anatomic geometry of
\textipa{[u]} that using four times as many elements does not change
the numerical result, and thus the resonances given in
Table~\ref{numberofelements} can be regarded as accurate in this
respect.

\begin{table} [t,h]
  \vspace{2mm}
  \centerline{
    \begin{tabular}{|c|c|c|c|c|c|}
      \hline
      Vowel & $R_1$ & $R_2$ & $R_3$ & $R_4$ & \# of elem.   \\
      \hline  
      \textipa{[\textscripta]} & 720 & 1547 & 2721 & 4138 & 47514  \\
      \textipa{[i]} & 246 & 2135 & 3592 & 4667 &  37335 \\
      \textipa{[u]} & 324 & 659 &  2262 & 3091 &  50579 \\
      \textipa{[\oe]} & 562 &  1612 &  2519 &  3602 & 53087 \\
      \hline
  \end{tabular}}
  \caption{\label{numberofelements} {\it Resonances (in Hz) of the
      Helmholtz problem \eqref{eq:eigenfuncs} by FEM, and the number
      of the elements in each geometry.}}
\end{table}

\section{Geometric data from MRI}

The raw MRI data was collected in pilot experiments in June 2010.  A
native male Finnish speaker pronounced prolonged vowels in a supine
position inside a MRI machine.  A data set of 53 simultaneously
recorded MRI data and sound samples was produced as reported in
\cite{A-H-M-P-P-S-V:RSSAMRI}.  For this study, four samples
corresponding to the vowels \textipa{[\textscripta]}, \textipa{[i]},
\textipa{[u]}, and \textipa{[\oe]} were chosen from this data set. The
selection criteria were 1) a visual quality assessment of the
spectrogram data, and 2) the requirement that the $F_1$-$F_2$ vowel
space should be covered in a satisfactory manner.  The spectrograms of
these samples can be found in \cite[p.64, p.66, p.68, and
  p.79]{Palo:L:2011}. The first three of these samples were pronounced
at the fundamental frequency $f_0 = 110$ Hz and the last one at $f_0 =
137.5$ Hz.

A Siemens Magnetom Avanto 1.5T scanner was used in these experiments.
A 12-element Head Matrix Coil was combined with a 4-element Neck
Matrix Coil in order to cover the anatomy of interest. 3D VIBE
(Volumetric Interpolated Breath-hold Examination)
\cite{R-L-L-P-K-T-A-W:AMRIVIBE} was found to be the most suitable MRI
sequence for rapid 3D acquisition. It was originally developed for
fast 3D imaging of the abdominal region where breath-hold during the
scan is essential, as the naming of the sequence suggests.  Sequence
parameters were optimized in order to minimize the acquisition time,
and we were able to carry out MRI with 1.8~mm isotropic voxels in just
7.6~s.

\begin{figure}[t]
\centerline{\includegraphics[width=50mm]{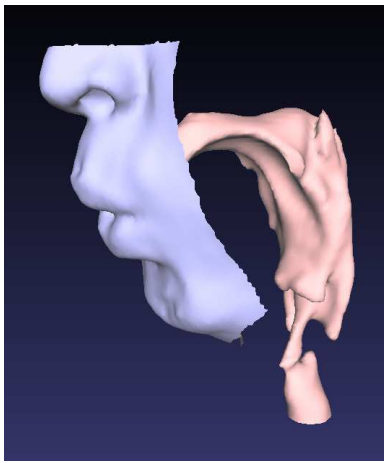}}
\caption{\label{fig:surface}{\it Air/tissue interface of the test subject while
    pronouncing \textipa{[y]}. The nasal cavities have been
    excluded.}}
\label{spprod}
\end{figure}

The tissue/air interface from the MR data was extracted by combining
sagittal DICOM sequence -images to form a 3D voxel image.  A
triangular surface mesh of the interface was then extracted using
custom MATLAB code.  The three boundary components $\Gamma_1$,
$\Gamma_2$, and $\Gamma_3$ were identified manually in the triangle
mesh so that the different boundary conditions could be applied in the
right places. Since teeth are not visible in MRI (and hence, they are
not part of the computational geometry of this work), some resulting
artefacts had to be corrected manually.  The velar port was open in
the geometries of \textipa{[\textscripta]} and \textipa{[i]}, and the
resulting hole in the surface model was manually closed.  A shaded
representation of a typical surface mesh is presented in
Fig.~\ref{fig:surface}.


The geometric error in the tissue/air interface is a fraction of the
voxel-based resolution of the original MRI data: interpolating in 2D
sections by the gray-scale values of pixels results in about 1.8 bits
of additional information compared to the MRI pixel size, corresponding
to the geometric error of $\approx$ 0.5 mm with the current
voxel size of 1.8 mm; see \cite{A-M-V-S-P:EMAEMRID}. We conclude that
the resonances in Table~\ref{numberofelements} do not contain
essential errors due to inaccuracies of surface geometries.

\section{Sound recording and formant extraction}

The interior of a MRI machine is a challenging environment for speech
recording.  We used the recording arrangement discussed in
\cite{M-P:RSDMRIII,Palo:L:2011} and the experimental arrangements
described in \cite[Section 2]{A-H-M-P-P-S-V:RSSAMRI}; see also
\cite{A-M-V-S-P:EMAEMRID}.  

Let us briefly describe the recording arrangement. A two-channel sound
collector samples the speech and noise signals in a dipole
configuration. The sound collector is an acoustically passive,
non-microphonic device which does not cause artifacts in the MR
images. The sound signals are coupled to a RF-shielded microphone
array by acoustic waveguides of length 3~m. There are tuned acoustic
impedance terminations at the both ends of the wave\-guides to
sufficiently reduce longitudinal resonances.  The microphone signals
are coupled to an amplifier that is situated outside the MRI room.
This analogue electronics is used to optimally subtract the noise
signal from the contaminated speech signal in real time, and the
cleaned-up signal is fed back to test subject's earphones.  The same
audio signal is digitized by a 24bit ADC, and the residual
longitudinal resonances of the wave\-guides are compensated
numerically in the post-processing stage.

\begin{table} [t,h]
  \vspace{2mm}
  \centerline{
    \begin{tabular}{|c|c|c|c|c|}
      \hline
      Vowel &  $F_1$ &  $F_2$ &  $F_3$ &  $F_4$  \\
      \hline  
      \textipa{[\textscripta]} 
      & 651 $\pm$ 7  &  1024 $\pm$ 35  &    2647 $\pm$ 117  &   3679 $\mp$ 36 \\
      \textipa{[i]} 
      & 247  $\pm$ 9 &  2183           & 3304 $\mp$ 46  &  4407 $\mp$ 251 \\
      \textipa{[u]} 
      & 306 $\mp$ 37  &   675 $\mp$ 39   &    2173 $\pm$ 13   &   3242 $\pm$ 139  \\
      \textipa{[\oe]}  
      & 483 $\mp$ 35 & 1249 $\pm$ 74 & 1994 $\mp$ 50 & 3188 $\mp$ 17  \\
      \hline
  \end{tabular}}
  \caption{\label{measuredformants} {\it Formants (in Hz) computed as
      means of those extracted from the beginnings and ends of the
      samples. The upper (lower) sign refers to the beginning (resp.,
      the end) of the sample.}}
\end{table}

For this work, the formants $F_1$...$F_4$ were first extracted
separately from the beginnings and the ends of the sound samples where
the acoustic MRI noise is absent. The extraction was done by LPC using
MATLAB similarly to the approach explained in \cite[Chapter 6]{Palo:L:2011}.
However, the present values were obtained by applying LPC analysis to FFT power
spectra of the signals with the algorithm detailed in \cite{Makhoul:SLPPA}. The
residual wave\-guide resonances were removed from the spectra before LPC
analysis. The results were compared visually to the peaks of the smoothed spectra to
detect crude errors. The final results in Table~\ref{measuredformants}
are the averages of these two values, augmented with their half
distances.  The formant data for \textipa{[i]} is subject to following
remarks: (i) the LPC found a peak at $855 \pm 133$ Hz but this was
removed from the data set as an outlier; (ii) $F_2$ could not be
extracted from the beginning sample by the LPC even though it can be
found easily in spectral curves by visual inspection; (iii) the LPC
finds very strong \emph{double} peaks about 500 Hz apart at $F_3$, and the
$F_3$-values given in Table~\ref{measuredformants} are defined as
their means.

\begin{table} [t,h]
  \vspace{2mm}
  \centerline{
    \begin{tabular}{|c|c|c|c|c|c|}
      \hline
      Vowel & $D_1$ & $D_2$ &  $D_3$ &  $D_4$ & mean discr.    \\
      \hline  
      \textipa{[\textscripta]} & 1.7  &  7.1 &  0.5  &   2.0 & 2.8 \\ 
      \textipa{[i]} & -0.1 &  -0.4   & 1.4  &  1.0 & 0.7\\
      \textipa{[u]} & 1.0 &  -0.4  &  0.7  & -0.8 & 0.7 \\
      \textipa{[\oe]} & 2.6 & 4.4 & 4.0 & 2.1 & 3.3 \\
      \hline
  \end{tabular}}
  \caption{\label{discrepancy} {\it Discrepancy (in semitones) between
      computed resonances and mean formant frequencies from
      Table~\ref{measuredformants}. Positive number implies that the
      computed resonance is higher than the measured formant.}}
\end{table}

The data of Table~\ref{measuredformants} can be found in
\cite{P-A-A-H-M-S-V:AFVRMRISD} except \textipa{[\oe]} (which is from a
different series where $f_0 = 137.5$ Hz) and $F_3$, $F_4$ of
\textipa{[i]} (which required a higher order LPC run).  These results
\emph{without} numerical compensation of the residual longitudinal
resonances of the waveguides can be found in \cite[Tables 6.2 and 6.3
  on p.49--50]{Palo:L:2011}.  We remark that a typical mean deviation
in vowel formant frequencies (when measured in ideal conditions in an
anechoic chamber rather than inside a MRI machine) is of the order of
0.5 semitones; see \cite{A-M-V-S-P:EMAEMRID}.

\section{Results and conclusions}

Just as in \cite{H-L-M-P:WFWE}, our new data indicates that the
computed resonances $R_i$ from \eqref{eq:eigenfuncs} tend to be higher
than the measured formants $F_i$.  The discrepancy given in
Table~\ref{discrepancy} is in semitone scale to make the comparison
easy with the ``$3\tfrac{1}{2}$ semitone rule'' that was discovered in
\cite{H-L-M-P:WFWE}. Recall that the difference of frequencies $F$ and
$R$ is $D = 12 \ln(R/F)/\ln(2)$ semitones.


The main sources of discrepancy in these computations and experiments
are as follows: (i) less than perfect performance of the test subject
in the MRI machine, (ii) sporadic problems in formant extraction by
the LPC, and (iii) physically unrealistic boundary conditions in
Eqs.~\eqref{eq:model}--\eqref{eq:eigenfuncs} especially at the mouth
opening.

The mean discrepancy in Table~\ref{discrepancy} is at its largest for
vowels \textipa{[\textscripta]} and \textipa{[\oe]} where the computed
resonances are consistently higher than the measured formants.  Also,
the mouth opening is largest for these vowels in our data set, and the
Dirichlet boundary condition on $\Gamma_1$ in \eqref{eq:eigenfuncs} is
expected to be the most significant error source. All this is in good
agreement with the results and the conclusions of \cite{H-L-M-P:WFWE}.

The particularly significant error in $F_2$ of
\textipa{[\textscripta]} can be explained by the fact that the velar port 
of the subject was unexpectedly open in the MR image, and the
anti-node of the standing pressure wave (corresponding to $F_2$) would
be at the velar port \emph{if} it was closed.  The nasality of the
pronunciation is clearly heard from the speech sample. In the
computational geometry, however, the hole was closed manually which
leads to the perfectly reflecting Neumann boundary condition for the
velar opening.  It is physically more realistic to use a similar
boundary condition for the velar opening as on $\Gamma_3$ in
\eqref{eq:eigenfuncs}.

It is worth noting that the computational model performs strikingly
well for \textipa{[u]}: The discrepancy is of the same order as the
fluctuations in formant values in sustained vowel production
\cite{A-M-V-S-P:EMAEMRID}.  The velar port is closed in this MRI geometry. Also,
the mouth opening is very small which results in relatively smaller
error due to the Dirichlet boundary condition.

We conclude that the error profile in Table \ref{discrepancy} supports
earlier observations in \cite{H-L-M-P:WFWE}, and it can be
qualitatively explained by considering the underlying physics.  A more
sophisticated exterior space model (compared to the Dirichlet boundary
condition) is likely to remove most of the formant discrepancy in
vowels where the mouth opening is large. Complications related to the
open velar port should be treated by taking into account the nasal tract
resonating structures. This can be done by including them in the
computational geometry or by setting an improved boundary condition at
the velar port opening.

The results of \cite{T-M-K:AAVTDVPFDTDM} support the view that
computed and measured resonances of a \emph{plastic model} VT do not
differ significantly from each other. Rather than carrying out speech
recordings in MRI, the authors produce 3D physical printouts from MRI
geometries by fast prototyping techniques for japanese vowels
\textipa{[a]}, \textipa{[i]}, \textipa{[u]}, \textipa{[e]},
\textipa{[o]}.  Separately imaged teeth geometries were manually
aligned with the soft tissue geometries.  The formant structure is
measured from the plastic models in ideal conditions, and the same
configuration is used for transfer function simulations by the
Finite-Difference Time-Domain (FDTD) method.

It is observed that the formant frequencies of the computed and
measured transfer functions in \cite{T-M-K:AAVTDVPFDTDM} differ from
each other by less than 3.2\% (i.e., 0.55 semitones) which is
comparable to our results on vowels \textipa{[i]} and
\textipa{[u]}. Such an indirect experimental arrangement excludes most
of the sources of discrepancy considered in our work, and the results
of \cite{T-M-K:AAVTDVPFDTDM} can therefore be regarded as a limit of
what is reasonable to expect in a computational modelling effort of a
natural vowel utterance.  The effect of anatomic details such as
piriform fossae, epiglottic valleculae and inter-dental spaces were
computed, and the contribution of the latter two was found to be
almost negligible in \cite{T-M-K:AAVTDVPFDTDM}.

\section{Acknowledgements}
The authors were supported by the Finnish Academy grant Lastu 135005,
the Finnish Academy projects 128204 and 125940, European Union grant
Simple4All, Aalto Starting Grant, and {\AA}bo Akademi Institute of
Mathematics.

 \newpage

\bibliographystyle{plain} 



\end{document}